\documentclass[12 pt,draft]{article}
\usepackage[cp1251]{inputenc}
\usepackage[russian]{babel}
\usepackage{amssymb,amsthm,amsmath}

\newtheorem{theorem}{Теорема}
\newtheorem*{corollary}{Следствие}

\newtheorem*{notation}{Замечание}
\newtheorem*{definition}{Определение}

\DeclareMathOperator{\Res}{Res} \DeclareMathOperator{\Log}{Log}
 
\DeclareMathOperator{\Int}{Int} 
\begin{document}

\noindent УДК 517.55
\\ \\
{\LARGE \textbf{Последовательности Риордана\\ и двумерные
разностные уравнения}} \footnote{Работа поддержана Российским
фондом фундаментальных исследований (РФФИ) и Государственным
фондом естественных наук Китая (ГФЕН) в рамках совместного проекта
<<Комплексный анализ и его приложения>> (проект N
08-01-92208\_ГФЕН)} {\LARGE\\ \\ \textbf{Riordan's arrays and\\
two-dimensional difference equations}}
\\ \\
\begin{flushright}
{\Large \textbf{А.П. Ляпин}}\footnote{e-mail: LyapinAP@yandex.ru}
\\
Институт педагогики, психологии и социологии\\
Сибирского федерального университета,\\
пр-т Свободный 79, Красноярск, 660041\\
\end{flushright}

{\noindent\small В работе предложено описание рациональных
последовательностей Риордана, возникающих в комбинаторном анализе,
как решений задачи Коши двумерных разностных уравнений
специального вида, и исследована асимптотика таких
последовательностей.

\noindent Ключевые слова: последовательности Риордана, многомерные
разностные уравнения, задача Коши, амеба характеристического
многочлена.
\\

\noindent In this paper the describing of rational Riordan's
arrays from combinatorial analysis is represented as solutions of
a Cauchy problem for two-dimensional difference equations and it
is researched the asymptotic of these arrays.

\noindent Keywords: Riordan's arrays, multidimensional difference
equations, Cauchy problem, amoeba of char\-act\-eri\-stic
polynomial.}

\section{Введение}

Понятие последовательности Риордана впервые появилось в работе
\cite{shapiro}  в связи с изучением групп Риордана и в дальнейшем
нашло широкое применение в таких задачах перечислительного
комбинаторного анализа, как задача о числе путей на целочисленной
решетке (\cite{merlini2003}), о производящих деревьях с
помеченными вершинами (<<level generating trees>>,
\cite{merlini2008}), о расстановке фигур на шахматной доске
(\cite{Egorychev}, \cite{Abramson}),  строки Блума (см.
\cite{merlini2008}, \cite{bloom1998}).

Приведем определение последовательности Риордана. Пусть $d(z) =
\sum_0^\infty d_k z^{-k-1}$ и $h(z) = \sum_0^\infty h_k z^{-k-1}$
~--- ряды Лорана, вообще говоря, формальные.
\begin{definition}
Последовательностью Риордана, ассоциированной с парой $d(z), h(z)$
называется последовательность $\{r(x,y), (x,y)\in
\mathbb{Z}_+^2\}$, производящая функция которой $$\mathcal{D}
(z,w) := \sum\limits_{(x,y)\in \mathbb{Z}_+^2}
\frac{r(x,y)}{z^{x+1} w^{y+1}}$$ имеет вид
\begin{equation}\label{riordan}
\mathcal{D} (z,w) = \frac {d(z)}{w-h(z)} := \sum_{y=0}^\infty
\frac{d(z) h^y(z)}{w^{y+1}}.
\end{equation}
\end{definition}
 Отметим, что $r(x,y) = \Res \{ d(z) h^y(z) z^x\}$, где $\Res$
~--- оператор, который ставит в соответствие некоторому
формальному ряду от $z$ коэффициент при $z^{-1}$.

Последовательности Риордана называются \emph{правильными}, если
$h_0\neq 0$ (см., например, \cite{merlini2003}).

Выделим класс последовательностей Риордана, которые будут
рассматриваться в данной работе.

\begin{definition}
Последовательность Риордана будем называются
\textbf{рациональной}, если ряд $h(z)$ является разложением в
окрестности бесконечно удаленной точки рациональной функции
$$
h(z) = \frac {Q(z)}{P(z)},
$$
где $P(z)=\sum\limits_{0}^{m} c_{\alpha 1}z^\alpha , Q(z) =
\sum\limits_\alpha c_{\alpha 0}z^\alpha$ ~--- многочлены от $z\in
\mathbb C$ и  $c_{m,1}\neq 0$
\end{definition}

Отметим, что из разложимости $h(z)$ в ряд Лорана в окрестности
бесконечно удаленной точки следует, что $\deg Q(z) < \deg P(z)=m$.

Рациональные последовательности Риордана будут \emph{правильными},
если справедливо равенство $\deg Q(z) +1 = \deg P(z)$
(эквивалентное условию $h_0 \neq 0$).

Основным результатом данной работы является описание рациональных
последовательностей Риордана как решения задачи Коши для одного
класса двумерных разностных уравнений. Подробная формулировка
задачи Коши, необходимые определения и обозначения приведены в п.
2.

\begin{theorem}\label{th1}
Двойная последовательность $\{r(x,y)\}$ является рациональной
последовательностью Риордана, определяемой парой $d(z)$ и $h(z) =
\frac {Q(z)}{P(z)}$, тогда и только тогда, когда она является
решением задачи Коши для разностного уравнения
\begin{equation}\label{razn_ur_th1}
\left[P(\delta_1) \cdot \delta_2 - Q(\delta_1)\right] r(x,y) = 0,
\end{equation}
\begin{equation}\label{nach_dan_th00}
r(x,y) = \varphi(x,y), (x,y)\ngeqslant (m,1),
\end{equation}
где функция $\varphi(x,y)$ задана на множестве $(x,y)\ngeqslant
(m,1)$ следующим образом:
\begin{equation}\label{nach_dan_th1}
\varphi(x,y) := \Res\left\{d(\xi) \left(\frac
{Q(\xi)}{P(\xi)}\right)^y \xi^x\right\}.
\end{equation}
\end{theorem}

Производящую функцию начальных данных \eqref{nach_dan_th00}
$$
\Phi (z,w) := \sum_{(x,y)\ngeqslant(m,1)}
\frac{\varphi(x,y)}{z^{x+1}w^{y+1}}
$$
после очевидной группировки всегда можно записать в виде
$$
\Phi(z,w) = \sum_{k=0}^{m-1} \frac {\Phi_k(w)}{z^{k+1}} + \frac
{d(z)}w, \text{ где } \Phi_k (w)  = \sum_{y=1}^{\infty}
\frac{\varphi(k,y)}{w^{y+1}},
$$
Следующая теорема позволяет представить производящую функцию
решения задачи Коши в виде простого выражения от $\Phi_k(w), d(z)$
и многочленов
$$
R_{k+1} (z,w) = \frac 1{z^{k+1}}\sum_{\alpha=k+1}^{m}
\left(c_{\alpha,1}
 w - c_{\alpha,0} \right)z^{\alpha},
$$
построенных на основе характеристического многочлена разностного
уравнения \eqref{razn_ur_th1}.

\begin{theorem}\label{th2}
Производящая функция решения задачи Коши
\eqref{razn_ur_th1}-\eqref{nach_dan_th00} имеет вид
$$
\mathcal D(z,w) = \left(P(z) d(z) + \sum_{k=0}^{m-1}
R_{k+1}(z,w)\Phi_k(w) - \frac 1w \sum_{\alpha=0}^{m-1}
\sum_{x=0}^{\alpha-1} \frac{c_{\alpha,0}
\varphi(x,0)}{z^{x-\alpha+1}} \right) / (P(z)w-Q(z)).
$$
\end{theorem}

\begin{corollary}
Необходимым и достаточным условием рациональности производящей
фун\-кции решения задачи Коши
\eqref{razn_ur_th1}-\eqref{nach_dan_th00} является рациональность
производящей функции $\Phi(z,w)$ начальных данных.
\end{corollary}

\begin{notation}
Если начальные данные определяют последовательность Риордана
(имеют вид \eqref{nach_dan_th1}), то в теореме \ref{th2} выражение
для производящей функции $\mathcal D(z,w)$ примет простой вид, так
как в этом случае имеет место равенство
$$
\sum_{k=0}^{m-1} R_{k+1}(z,w)\Phi_k(w) - \frac 1w
\sum_{\alpha=0}^{m-1} \sum_{x=0}^{\alpha-1} \frac{c_{\alpha 0}
\varphi(x,0)}{z^{x-\alpha+1}} =0.
$$
\end{notation}

Следующая теорема позволяет найти асимптотику последовательностей
Риордана. Вопрос об отыскании асимптотического поведения решений
многомерных разностных уравнений весьма актуален и рассматривался
многими авторами в работах \cite{Egorychev}, \cite{tz91},
\cite{orlov}, \cite{pem}, \cite {wilson}, \cite{LPZ2005},
\cite{LPZ2008}.

Одним из способов исследования свойств двумерной
последовательности $r(x,y)$ является изучение асимптотики ее
<<диагональных>> подпоcледовательностей. А именно, для двойных
последовательностей $\{r(x,y)\}$ определим диагональные
подпоследовательности следующим образом: фиксируем $(p,q)\in
\mathbb{Z}_{+}^{2}$ и будем рассматривать одномерную
последовательность
$\{r(p\lambda,q\lambda)\}_{\lambda\in\mathbb{Z}_{+}}$.

Заметим, что такой подход к изучению асимптотического поведения
двойной последовательности применялся в \cite{tz91} при решении
проблемы устойчивости двумерных цифровых рекурсивных фильтров. В
работе \cite{pem} асимптотика коэффициентов рациональной
производящей функции изучалась в связи с такими задачами
перечислительного комбинаторного анализа, как, например, задача о
подбрасывании несимметричной монеты.

В случае рациональных последовательностей Риордана условия,
обеспечивающие применимость метода перевала удобно сформулировать
в терминах многогранника Ньютона $\mathcal N_R$ и амебы
$\mathcal{A}_R$ характеристического многочлена $R$, подробные
определения и обозначения которых приведены в п. 2.

Отметим лишь, что $E_{m,1}$ --- это компонента дополнения амебы,
соответствующая вершине многогранника Ньютона $(m,1)$, для которой
сформулирована задача Коши, а через $\Omega_{m,1}$ обозначен
конус, порожденный векторами $\{(m,1)-(\tau_1, \tau_2)\}$, где
$(\tau_1, \tau_2)\in \mathcal N_R$.

\begin{theorem}\label{th3}
Пусть функция $d(z)$ голоморфна вне особых точек функции $h(z)=
\frac QP$. Если все корни многочленов $P(z)$ и $Q(z)$ простые и
различные по модулю (каждого многочлена в отдельности) и граница
амебы характеристического многочлена $R(z,w) = P(z)\cdot w-Q(z)$
гладкая, тогда для всякого направления $(p,q)\in\Int \Omega_{m,1}$
существует единственная точка $(z_0(\frac pq),w_0(\frac pq))$,
такая, что ее логарифмический образ $\Log(z_0,w_0)\in
\partial E_{m,1}$, и
\begin{align}\label{assymp}
r(x,y) \sim \frac {d(z_0)}{\sqrt{2\pi \lambda q H(z_0)}} \left[
z_0^{p}w_0^{q}\right]^\lambda, x = \lambda p, y = \lambda
q,\lambda \to \infty,
\end{align}
где
$$H(z) = \frac{Q''(z)}{Q(z)} -\frac{P''(z)}{P(z)} + 2 \frac pq \frac 1z
\frac{P'(z)}{P(z)} - \frac pq (1+\frac pq) \frac 1{z^2}.$$
\end{theorem}

Отметим, что точка $(z_0,w_0)$ удовлетворяет системе уравнений
\begin{align*}
\begin{cases}
P(z) w - Q(z) = 0\\
z\left(\frac{P'(z)}{P(z)} - \frac{Q'(z)}{Q(z)} \right) = \frac pq
\end{cases}.
\end{align*}

\section{Определения и обозначения}

Пусть $\mathbb{Z}^{n}$ --- $n$-мерная целочисленная решетка, точки
которой обозначим $x = (x_1, \ldots, x_n)$, $\alpha = (\alpha_1,
..., \alpha_n)$, $\mathbb{Z}^{n}_+$ ~--- ее положительный октант,
множество $C \subset \mathbb Z^n_+$ конечное. Определим на
комплекснозначных функциях $r: \mathbb Z^{n}\to \mathbb{C}$
операторы сдвига $\delta_j$ вида
\begin{align*}
\delta_j r(x) = r(x_1, \ldots, x_j+1, \ldots, x_n),\, j=1, \ldots,
n,
\end{align*}
и рассмотрим полиномиальный разностный оператор с постоянными
коэффициентами
$$
R(\delta) = R(\delta_1, \delta_2, ..., \delta_n) = \sum_{\alpha\in
C} c_\alpha \delta^\alpha = \sum_{(\alpha_1,...,\alpha_n)\in C}
c_{\alpha_1,...,\alpha_n} \delta_1^{\alpha_1} \cdots
\delta_n^{\alpha_n}
$$

Нас интересуют (однородные) разностные уравнения вида $R(\delta)
r(x) = 0$ или
\begin{align} \label{razn_ur}
\sum_{\alpha\in C} c_\alpha r(x+\alpha) = 0,
\end{align}
 где $r(x)$ --- неизвестная функция целочисленных аргументов
$x = (x_1, \ldots, x_n)$.

Пусть $l=(l_1, ..., l_n) \in \mathbb Z^n_+, l \neq 0$ ~---
некоторая точка $n$-мерной целочисленной решетки $\mathbb Z^n_+$ и
на множестве $X_l = \{ x\in \mathbb Z ^n _+ \setminus (l+\mathbb Z
^n _+)\}$ задана функция
\begin{align}\label{nach_dan}
\varphi: X_l \to \mathbb C.
\end{align}

В дальнейшем будем  вместо $x\in \mathbb Z ^n _+ \setminus
(l+\mathbb Z ^n _+)$ использовать запись $x \ngeqslant l$.

\textbf{Задача Коши} состоит в отыскании решения разностного
уравнения \eqref{razn_ur}, совпадающего на множестве $X_l$ с
функцией \eqref{nach_dan}.

Условия на множество $X_l$, обеспечивающие существование и
единственность решения задачи Коши приведены в \cite{petrovsek}, в
\cite{lein2007} эти условия сформулированы в терминах
многогранника Ньютона характеристического многочлена $R$.

Характеристическим многочленом разностного уравнения
\eqref{razn_ur} называется многочлен вида $ R(z)= \sum\limits
_{\alpha\in C} c_\alpha z^\alpha.$

\emph{Многогранник Ньютона} $\mathcal N_R$ характеристического
многочлена $R(z)$ представляет выпуклую оболочку в
$\mathbb{R}^{n}$ конечного множества точек $C$.

Решение задачи Коши \eqref{razn_ur_th1}-\eqref{nach_dan_th00}
существует и единственно, если для некоторой (целочисленной) точки
$m= (m_1, ..., m_n)\in \mathcal N_R$ выполняется $\mathcal N_R
\subseteq \Pi_m$, где $\Pi_m=\{\alpha\in\mathbb R^n_+:
\alpha_j\leqslant m_j, j=1,...,n\}$.

\emph{Амебой $\mathcal{A}_R$ многочлена $R(z)$} называется образ
множества нулей этого многочлена $$ \mathcal{V} = \{z\in
\mathbb{C}^{n} : R(z)=0\}$$ при логарифмическом проектировании
$\Log: (z_1,\ldots,z_{n}) \mapsto (\log |z_1|,\ldots,\log
|z_{n}|).$

Дополнение к амебе состоит из конечного числа связных компонент,
ограниченного снизу числом вершин многогранника Ньютона, а сверху
--- числом целых точек пересечения $\mathcal{N}_R \cap
\mathbb{Z}^{n}$. Кроме того (\cite{paasare}), каждой вершине $\nu$
многогранника Ньютона можно сопоставить непустую связную
компоненту $E_\nu$ дополнения амебы $\mathcal{A}_R$ и разложение в
ряд Лорана функции $1/R(z)$, сходящееся в $\Log^{-1}E_\nu$.

\section{Доказательства}

\begin{proof}[Доказательство теоремы 1] Докажем необходимое условие.
Сгруппируем ряд $\mathcal D(z,w)$ по отрицательным степеням $w$:
\begin{align*}
 \mathcal D (z,w) =
 \frac{d(z)}{w\left(1-\frac{Q(z)}{wP(z)}\right)} =
 \sum_{y=0}^\infty d(z) \left(\frac {Q(z)}{P(z)}\right)^y \frac
 1{w^{y+1}},
\end{align*}
откуда следует, что $r(x,y) = \Res \left\{ d(z)
\left(\frac{Q(z)}{P(z)}\right)^y z^x\right\}$ (см. замечание к
определению последовательности Риордана).

В тождестве
$$[P(z)\cdot w - Q(z)] \mathcal D(z,w) \equiv P(z) d(z)$$
преобразуем левую часть:
\begin{align*}
&\left( P(z)\cdot w - Q(z)\right) \cdot \sum_{x,y\geqslant 0}
\frac{r(x,y)}{z^{x+1}w^{y+1}} = \\
&= P(z) w \left( \sum_{x\geqslant 0,y=0}
\frac{r(x,y)}{z^{x+1}w^{y+1}} + \sum_{x\geqslant 0, y \geqslant 1}
\frac{r(x,y)}{z^{x+1}w^{y+1}}\right) - Q(z) \sum_{x,y\geqslant 0}
\frac{r(x,y)}{z^{x+1}w^{y+1}} = \\ &= \sum_{\alpha=0}^m
c_{\alpha,1} z^\alpha \sum_{x=0}^\infty
\frac{\varphi(x,0)}{z^{x+1}} + \sum_{\alpha=0}^{m} c_{\alpha,1}
z^\alpha \left[\sum_{\substack{x=0, \ldots, \alpha-1\\y\geqslant
0}} \frac{\varphi(x,y+1)}{z^{x+1}w^{y+1}}+
\sum_{\substack{x\geqslant \alpha\\y\geqslant 0}}
\frac{r(x,y+1)}{z^{x+1}w^{y+1}} \right] -\\&- \sum_{\beta=0}^{m}
c_{\beta,0} z^\beta \left[ \sum_{\substack{x=0, \ldots,
\beta-1\\y\geqslant 0}} \frac{\varphi(x,y)}{z^{x+1}w^{y+1}} +
\sum_{\substack{x\geqslant \beta\\ y\geqslant 0}}
\frac{r(x,y)}{z^{x+1}w^{y+1}} \right] = \\ &= P(z)d(z) +
\sum_{\alpha=0}^{m} \left(c_{\alpha,1}z^\alpha
\sum_{\substack{x=0, \ldots, \alpha-1\\y\geqslant 0}}
\frac{\varphi(x,y+1)}{z^{x+1}w^{y+1}}\right) -
\sum_{\beta=0}^{m}\left( c_{\beta,0}z^\beta \sum_{\substack{x=0,
\ldots, \beta-1\\y\geqslant 0}} \frac{\varphi(x,y)}{z^{x+1}w^{y+1}}\right) + \\
&+ \sum_{\substack{x\geqslant 0\\y\geqslant 0}} \left[
\sum_{\alpha=0}^{m} c_{\alpha,1} r(x+\alpha, y+1) -
\sum_{\beta=0}^{m} c_{\beta,0} r(x+\beta,y)\right] \frac1{z^{x+1}
w^{y+1}}.
\end{align*}

Ввиду отсутствия в правой части тождества переменных в
отрицательной степени, получаем для всех $x\geqslant 0, y\geqslant
0$ соотношение
\begin{align*}
\sum_{\substack{x\geqslant 0\\y\geqslant 0}} \left[
P(\delta_1)\cdot \delta_2  - Q(\delta_1) \right] r(x,y) =0,
\end{align*}
что и завершает доказательство необходимого условия теоремы.

Докажем достаточность. Для доказательства нам потребуется
следующее утверждение:
\begin{align}\label{utv9}
\sum_{(\alpha, \beta)\in A} c_{\alpha,\beta} z^\alpha w^\beta
\sum_{(x,y)\ngeqslant (\alpha, \beta)} \frac{\xi^x
\eta^y}{z^{x+1}w^{y+1}} = \frac 1{(z-\xi)(w-\eta)}\sum_{(\alpha,
\beta)\in A} c_{\alpha,\beta} (z^\alpha w^\beta - \xi^\alpha
\eta^\beta).
\end{align}

Запишем разностное уравнение \eqref{razn_ur_th1} в виде
$$\sum_{\alpha=0}^{m} c_{\alpha,1} r(x+\alpha, y+1) -
\sum_{\alpha=0}^{m} c_{\alpha,0} r(x+\alpha,y) =0,$$ домножим его
левую часть на моном $z^{-x-1}w^{-y-1}$ и просуммируем по всем
$x\geqslant 0, y\geqslant 0$:
\begin{align*}
&\sum_{\substack{x\geqslant 0\\y\geqslant 0}}
\left(\sum_{\alpha=0}^{m} c_{\alpha,1} r(x+\alpha, y+1) -
\sum_{\alpha=0}^{m} c_{\alpha,0}
r(x+\alpha,y)\right)\frac{1}{z^{x+1}w^{y+1}}=\\
=&\sum_{\substack{x\geqslant 0\\y\geqslant 0}} \sum_{\alpha=0}^{m}
\frac{c_{\alpha,1} r(x+\alpha, y+1)z^\alpha
w}{z^{x+\alpha+1}w^{y+2}} - \sum_{\substack{x\geqslant
0\\y\geqslant 0}} \sum_{\alpha=0}^{m} \frac{c_{\alpha,0}
r(x+\alpha, y)z^\alpha
}{z^{x+\alpha+1}w^{y+1}}=\\=&\sum_{\alpha=0}^{m}c_{\alpha,1}z^\alpha
w \sum_{\substack{x\geqslant \alpha\\y\geqslant 1}}
\frac{r(x,y)}{z^{x+1}w^{y+1}} -
\sum_{\alpha=0}^{m}c_{\alpha,0}z^\alpha \sum_{\substack{x\geqslant
\alpha\\y\geqslant 0}} \frac{r(x,y)}{z^{x+1}w^{y+1}}=\\=&
\sum_{\alpha=0}^{m}c_{\alpha,1}z^\alpha w \left(\mathcal D(z,w) -
\sum_{(x,y)\ngeqslant(\alpha,1)}
\frac{r(x,y)}{z^{x+1}w^{y+1}}\right) -
\sum_{\alpha=0}^{m}c_{\alpha,0}z^\alpha \left(\mathcal D(z,w) -
\sum_{(x,y)\ngeqslant(\alpha,0)}
\frac{r(x,y)}{z^{x+1}w^{y+1}}\right) = \\=&\mathcal D(z,w) R(z,w)
- \left(\sum_{\alpha=0}^{m}c_{\alpha,1}z^\alpha w
\sum_{(x,y)\ngeqslant(\alpha,1)} \frac{r(x,y)}{z^{x+1}w^{y+1}} -
\sum_{\alpha=0}^{m}c_{\alpha,0}z^\alpha
\sum_{(x,y)\ngeqslant(\alpha,0)}
\frac{r(x,y)}{z^{x+1}w^{y+1}}\right) \tag{*}
\end{align*}
Поскольку на множестве начальных данных функция $\varphi(x,y)$
имеет вид \eqref{nach_dan_th1}, то в силу линейности оператора
$\Res$ (см. \cite[стр. 15]{Egorychev}) и формулы \eqref{utv9}:
\begin{align*}
&\sum_{\alpha=0}^m c_{\alpha,1}z^\alpha w
\sum_{(x,y)\ngeqslant(\alpha,1)} \frac{r(x,y)}{z^{x+1}w^{y+1}}
=\\=& \Res \left\{\frac{d(\xi) P(\xi)}{P(\xi) \eta -
Q(\xi)}\frac{1}{(w-\eta)(z-\xi)}\left(\sum_{\alpha=0}^m
c_{\alpha,1}z^\alpha w - \sum_{\alpha=0}^m c_{\alpha,1}\xi^\alpha
\eta\right)\right\} = \\ =& \Res \left\{\frac{d(\xi)
P(\xi)}{P(\xi) \eta - Q(\xi)}\frac{P(z)w -
P(\xi)\eta}{(w-\eta)(z-\xi)}\right\},
\end{align*}
аналогично,
\begin{align*}
&\sum_{\alpha=0}^m c_{\alpha,1}z^\alpha w
\sum_{(x,y)\ngeqslant(\alpha,1)} \frac{r(x,y)}{z^{x+1}w^{y+1}} =
\Res \left\{\frac{d(\xi) P(\xi)}{P(\xi) \eta - Q(\xi)}\frac{Q(z) -
Q(\xi)}{(w-\eta)(z-\xi)}\right\}.
\end{align*}
Следовательно, в силу свойств оператора $\Res$,
\begin{align*}
&\sum_{\alpha=0}^{m}c_{\alpha,1}z^\alpha w
\sum_{(x,y)\ngeqslant(\alpha,1)} \frac{r(x,y)}{z^{x+1}w^{y+1}} -
\sum_{\alpha=0}^{m}c_{\alpha,0}z^\alpha
\sum_{(x,y)\ngeqslant(\alpha,0)} \frac{r(x,y)}{z^{x+1}w^{y+1}}
=\\=& \Res \left\{ \frac{d(\xi) P(\xi)}{P(\xi) \eta -
Q(\xi)}\frac{P(z)w - P(\xi)\eta}{(w-\eta)(z-\xi)} - \frac{d(\xi)
P(\xi)}{P(\xi) \eta - Q(\xi)}\frac{Q(z) -
Q(\xi)}{(w-\eta)(z-\xi)}\right\}=\\=& \Res \left\{ \frac{d(\xi)
P(\xi)}{P(\xi) \eta - Q(\xi)}\frac {(P(z)w - Q(z)) - (P(\xi)\eta -
Q(\xi))}{(w-\eta)(z-\xi)}\right\} =\\=& (P(z)w - Q(z))
\Res\left\{\frac{d(\xi) P(\xi)}{(P(\xi) \eta -
Q(\xi))(w-\eta)(z-\xi)}\right\} - \Res \left\{\frac{d(\xi)
P(\xi)}{(w-\eta)(z-\xi)}\right\} = P(z) d(z).
\end{align*}

Используя это равенство, получим  $$ \mathcal D(z,w) R(z,w) - P(z)
d(z) = 0, $$ откуда сразу следует, что $r(x,y)$ --- рациональная
последовательность Риордана.
\end{proof}

\begin{proof}[Доказательство теоремы 2]
В соотношении (*) из  доказательства теоремы 1
\begin{align*}
&\mathcal D(z,w)R(z,w) =\\=&
\sum_{\alpha=0}^{m}c_{\alpha,1}z^\alpha w
\sum_{(x,y)\ngeqslant(\alpha,1)} \frac{r(x,y)}{z^{x+1}w^{y+1}} -
\sum_{\alpha=0}^{m}c_{\alpha,0}z^\alpha
\sum_{(x,y)\ngeqslant(\alpha,0)}
\frac{r(x,y)}{z^{x+1}w^{y+1}}=\\=&
\sum_{\alpha=0}^{m}c_{\alpha,1}z^\alpha w \left( \sum_{x=0}^\infty
\frac{\varphi(x,0)}{z^{x+1}w} +
\sum_{x=0}^{\alpha-1}\sum_{y=1}^\infty
\frac{\varphi(x,y)}{z^{x+1}w^{y+1}}\right) -\\-&
\sum_{\alpha=0}^{m}c_{\alpha,0}z^\alpha \left(
\sum_{x=0}^{\alpha-1} \frac{\varphi(x,0)}{z^{x+1}} +
\sum_{x=0}^{\alpha-1}\sum_{y=1}^\infty
\frac{\varphi(x,y)}{z^{x+1}w^{y+1}}\right) =\\=& P(z)d(z) +
\sum_{\alpha=0}^m c_{\alpha,1}z^\alpha w
\sum_{x=0}^{\alpha-1}\frac{\Phi_x(w)}{z^{x+1}} - \sum_{\alpha=0}^m
c_{\alpha,0}z^\alpha
\sum_{x=0}^{\alpha-1}\frac{\varphi(x,0)}{z^{x+1}w} -
\sum_{\alpha=0}^m c_{\alpha,0}z^\alpha
\sum_{x=0}^{\alpha-1}\frac{\Phi_x(w)}{z^{x+1}}=\\=& P(z)d(z) +
\sum_{\alpha=0}^m \sum_{x=0}^{\alpha-1}
(c_{\alpha,1}w-c_{\alpha,0})z^\alpha
\frac{\Phi_x(w)}{z^{x+1}}-\frac 1w\sum_{\alpha=0}^m
\sum_{x=0}^{\alpha-1}
\frac{c_{\alpha,0}\varphi(x,0)}{z^{x-\alpha+1}}=\\=& P(z)d(z) +
\sum_{x=0}^{m-1} \left( \Phi_x(w) \cdot \frac 1{z^{x+1}}
\sum_{\alpha=x+1}^{m}(c_{\alpha,1}w-c_{\alpha,0})z^\alpha\right)-\frac
1w\sum_{\alpha=0}^m \sum_{x=0}^{\alpha-1}
\frac{c_{\alpha,0}\varphi(x,0)}{z^{x-\alpha+1}} =\\=&P(z)d(z) +
\sum_{x=0}^{m-1} \Phi_x(w) \cdot R_{x+1}(z,w) - \frac
1w\sum_{\alpha=0}^m \sum_{x=0}^{\alpha-1}
\frac{c_{\alpha,0}\varphi(x,0)}{z^{x-\alpha+1}}.
\end{align*}
Откуда и следует утверждение теоремы.
\end{proof}

\begin{proof}[Доказательство теоремы \ref{th3}]
Пусть $a_1, a_2, ..., a_{N_1}$ --- корни многочлена $Q(z)$, $b_1,
b_2, ..., b_{N_2}$ --- корни многочлена $P(z)$. Так как согласно
условиям теоремы, все корни простые и различные по модулю, тогда
амебу  $\mathcal A_R$ характеристического многочлена можно
представить как множество точек $\mathbb{R}^2_{(\xi,\eta)}$ вида
$\xi=t, \eta = f(t, \varphi)$, где
\begin{align*}
f(t,\varphi)=\log \frac {|e^{t+i\varphi}-a_1|\cdot \ldots \cdot
|e^{t+i\varphi}-a_{N_1}|}{|e^{t+i\varphi}-b_1|\cdot \ldots \cdot
|e^{t+i\varphi}-b_{N_2}|}, (t,\varphi)\in\mathbb R\times\left[0,
2\pi\right).
\end{align*}

Отметим, что каждому $a_i\neq 0$ и $b_j\neq 0$ соответствуют лучи
на прямых $\xi = \log |a_i|$ и $\xi = \log |b_j|$, целиком
принадлежащие амебе. Также амеба содержит еще два луча,
возникающих при $t\to\pm\infty$. Таким образом, число компонент
дополнения к амёбе не меньше $N_1+N_2+2-\kappa$, где $N_1$ и $N_2$
--- число различных по модулю и отличных от ноля корней
многочленов $Q(z)$ и $P(z)$ соответственно, а $\kappa=1$, если у
одного из многочленов есть корень, равный нолю, и $\kappa=0$ в
остальных случаях. С другой стороны, так как все целочисленные
точки многогранника Ньютона $\mathcal{N}_R$ характеристического
многочлена $R$ лежат на его границе, то их число равно
$N_1+N_2+2-\kappa$.

Таким образом, число компонент дополнения амебы равно числу целых
точек в многограннике Ньютона, это и означает по определению, что
амеба характеристического многочлена максимальна.

Условие гладкости границы позволит применить теорему 2 из работы
\cite{LPZ2008}, из которой для $(p,q)\in\Int \Omega_{m,1}$ следует
справедливость формулы $$r(x,y) \sim \frac {C(p,q)}{\sqrt{2\pi
\lambda}} \left[ z_0^{p}w_0^{q}\right]^\lambda, x = \lambda p, y =
\lambda q,\lambda \to \infty,$$ где константа вычисляется по
формуле $$ C(p,q) = \frac {d(z_0)}{\sqrt{q H(z_0)}}, H(z) =
\frac{Q''(z)}{Q(z)} -\frac{P''(z)}{P(z)} + 2 \frac pq \frac 1z
\frac{P'(z)}{P(z)} - \frac pq (1+\frac pq) \frac 1{z^2},$$ а точка
$(z_0,w_0)$ удовлетворяет системе
\begin{align*}
\begin{cases}
P(z) w - Q(z) = 0\\
z\left(\frac{P'(z)}{P(z)} - \frac{Q'(z)}{Q(z)} \right) = \frac pq
\end{cases}.
\end{align*}
\end{proof}

\section{Примеры}

Приведем несколько примеров рациональных последовательностей
Риордана.

\textbf{Пример 1.} \emph{Биномиальные коэффициенты} являются
решением разностного уравнения
\begin{align*}
    f(x+1, y+1) - f(x,y+1) - f(x,y)  =0,
\end{align*}
с начальными данными
\begin{equation*}
   \varphi (x,y) =
   \begin{cases}
        1,     \text{если} &x \geqslant  0, y =0\\
        0,     \text{если} &x=0, y \geqslant 1
   \end{cases}.
\end{equation*}
Производящая функция двойной последовательности $\{f(x,y)\}$ имеет
вид:
$$
\mathcal D(z,w) = \frac 1{zw-w-1} = \frac {d(z)} {w - h(z)},
$$
где
$$
d(z) = h(z) = \frac 1{z-1}.
$$
Тогда по теореме \ref{th3} для направлений $(p,q), \frac pq > 1$
справедлива асимптотическая формула:
$$
r(x,y) \sim \frac 1q \sqrt{\frac {pq}{2\pi \lambda (p-q)}}
\left((p-q)^{p-q} \frac {p^p}{q^q}\right)
$$

\textbf{Пример 2.} Пусть $A=(a_{ij})$ --- это $x\times m$ матрица
с $x\cdot m$ различными элементами. Обозначим $r_m (x,y),
0<y\leqslant x$ --- число способов выбора $y$ из $x\cdot m$
элементов таким образом, чтобы никакие два не стояли в одной
строке, а если выбраны элементы из смежных строк, то они должны
стоять в одном столбце.

В работе \cite{Abramson} найдена явная формула для вычисления $r_m
(x,y)$:
$$
r_m (x,y) = \sum_{r=1}^{y} m^r \binom{y-1}{r-1} \binom{x-y+1}{r}
$$
и приведено разностное уравнение для $r_m(x,y)$:
$$
r_m (x+2,y+1) - r_m (x+1,y+1)-r_m (x+1,y)-(m-1)r_m (x,y) =0.
$$
Рассмотрим задачу Коши для данного уравнения, указав <<начальные>>
значения на множестве $X_{2,1} = \{(x,y)\in\mathbb{Z}^2_+ : (x,y)
\ngeqslant (2,1) \}$:
\begin{equation*}
   \varphi (x,y) =
   \begin{cases}
        1,   \text{если} &x\geqslant 0, y=0,\\
        0,   \text{если} &x=0, y \geqslant 1 \text{ или } x=1, y\geqslant 2\\
        m,   \text{если} &(x,y)=(1,1)
   \end{cases},
\end{equation*}

По теореме \ref{th2} найдем производящую функцию для коэффициентов
$r_m (x,y)$:
\begin{align*}
\mathcal D(z,w)= \frac{P(z) d(z)}{P(z)w-Q(z)} =
\frac{z}{z(z-1)w-z-(m-1)}.
\end{align*}

Используя теорему 3, вычислим асимптотику $r(x,y)$ при $m=2$.

Точка $(z_0,w_0)$ является решением системы
\begin{equation*}
   \begin{cases}
        w=\frac{z+1}{z(z-1)}\\
        z\left(\frac{2z-1}{z(z-1)}-\frac1{z+1}\right)=\mu
   \end{cases}.
\end{equation*}
Наибольший вклад в асимптотику даст точка $(z_0,w_0)$ вида
$$\left(\frac{1+
\sqrt{(\mu-1)^2+1}}{\mu-1},
\frac{(\mu-1)(\sqrt{2-2\mu+\mu^2}+\mu)}
    {(2+\sqrt{2-2\mu+\mu^2}-\mu)(1+\sqrt{2-2\mu+\mu^2})}\right) $$
Далее, $$H(z_0) = \frac{(\mu-1)((\mu-2)z_0-\mu)}{z_0^2(z_0-1)},
d(z_0)=\frac 1{z_0-1}$$ и, применяя формулу \eqref{assymp}, для
$\frac pq>1$ получим $$ r(x,y) \sim \frac{z_0}{\sqrt{2\pi \lambda
q (\mu-1)((\mu-2)z_0-\mu)(z_0-1)}}\left[z_0^p
w_0^q\right]^\lambda, x = \lambda p, y=\lambda q, \lambda \to
\infty.$$

\textbf{Пример 3.} Рассмотрим последовательности из $n$ элементов
$a_1 a_2 ... a_n$, причем $a_1=0$ и $a_j\in\{0,1\}$ для
$2\leqslant j \leqslant n$. Элемент последовательности $a_j$
назовем изолированным, если он отличен от всех элементов, стоящих
на соседних местах. Обозначим $r(n,k)$ число таких
последовательностей, содержащих $k$ изолированных элементов.
Очевидно, что $r(n,k) =0$, если $n<k$ (см. \cite{bloom1998},
\cite{merlini2008}).

Данная последовательность является решением задачи Коши для
разностного уравнения
$$
r(x+2, y+1) - r(x+1, y+1) - r(x+1, y) - r(x, y+1) + r(x,y) =0
$$
с начальными данными $\varphi(0,0) = 1, \varphi(1,0) = 0,
\varphi(x,0) = \varphi(x-1,0)+\varphi(x-2,0), x\geqslant 2$,
$\varphi(1,1)=1$, $\varphi(0,y) = 0, y\geqslant 1$ и $\varphi(1,y)
= 0, y\geqslant 2$.

По теореме \ref{th2} производящая функция имеет вид
$$
\mathcal D (z,w) = \frac {z - 1}{z^2w - zw - w - z + 1} = \frac
{d(z)}{w-h(z)},
$$
где
$$
d(z) = h(z) = \frac {z-1}{z^2-z-1}.
$$

\end{document}